\documentclass[doublespacing]{elsart}
\usepackage{amssymb}
\usepackage[english,francais]{babel}
\newtheorem{e-proposition}[theorem]{Proposition}

\newtheorem{e-definition}[theorem]{Definition\rm}

\newtheorem{theoreme}{Th\'eor\`eme}[section]

\newtheorem{proposition}[theoreme]{Proposition}

\newcommand{\R}{I\!\! R}%For symbol of Euclidean space R

\def\og{\leavevmode\raise.3ex\hbox{$\scriptscriptstyle\langle\!\langle$~}}
\def\fg{\leavevmode\raise.3ex\hbox{~$\!\scriptscriptstyle\,\rangle\!\rangle$}}
\begin{document}
\centerline{}
\begin{frontmatter}
\selectlanguage{english}
\title{ ON COMPACT FINSLER SPACES OF POSITIVE CONSTANT CURVATURE }
\selectlanguage{english}
\author{Behroz Bidabad},
\ead{bidabad@aut.ac.ir}
\address{Department of Mathematics, Amirkabir University of Technology \\(Tehran Polytechnic), Tehran 15914, Iran.}
%\medskip
%\begin{center}
%{\small %Received  ...2011; accepted after revision ....\\
%Presented by Marcel Berger}
%\end{center}
%-------------------------------------------------------------------------------------------
\begin{abstract}
An $n$-dimensional ($n\geq 2$)  simply connected, compact without boundary Finsler space  of positive constant sectional curvature is  conformally homeomorphic to  an n-sphere in the Euclidean space $\R^{n+1}$.
 \vskip
0.5\baselineskip
\noindent{\bf R\'esum\'e} \vskip 0.5\baselineskip \noindent {\bf Sur les espaces   finsl\'eriennes  compactes    \`a courbure  constante  positive.}
\selectlanguage{francais}
 Un espace de Finsler de dimension $n$ ($n\geq 2$), simplement connexe  compacte, non-born\'ee,   \`a courbure  sectionelle constante  positive
 est conform\'ement hom\'eomorphe \`a une $n$-sph\`ere  dans l'espace euclidienne
 $\mathbb{R}^{n+1}$.
\end{abstract}
%-------------------------------------------------------------------------------------------
\end{frontmatter}
\selectlanguage{english}
{\small\emph{ Mathematics Subject Classification}: Primary 53C60; Secondary 58B20.\\
%\emph{ Subject Classification}: Finsler manifolds, Real and complex differential geometry, Global analysis.}\\
\section*{Introduction.}
The normal coordinates has been proved to be an extremely useful tool in the theory of global
Riemannian
geometry, while it is not so useful in Finsler geometry. In fact, in
the latter case the exponential map is only $C^1$ at the zero
section of $TM$ while in the former case it is $C^\infty$. More intuitively H. Akbar-Zadeh proved that exponential map is $C^2$ at the zero
section if and only if the Finsler structure is of Berwald type, cf. \cite{AZ}. In 1950s  Y. Tashiro have defined  adapted coordinates which was used by many authors to obtain global results in Riemannian geometry, which will not be mentioned here. Effectively, definition of this coordinate system on a manifold $M$, is equivalent to the existence of a non-trivial solution for a certain second order differential equation which describes circles preserving transformations on  $M$. Recently, the circle-preserving transformations are studied in Finsler geometry by present author and Z. Shen, cf. \cite{BS}.
Previously, inspiring the work of Tashiro, cf. \cite{T}, the present author have specialized adapted coordinates to Finsler setting in a joint paper, and proved the following theorem, cf. \cite{AB}. \\
\textit{\textbf{Theorem A:}}\emph{Let $(M,g)$ be a connected complete
Finsler manifold of dimension $n\geq 2$. If $M$ admits a non-trivial
solution of
\begin{equation}
 \label{c-field}
\nabla_i\nabla_j\rho= \phi g_{ij},
 \end{equation}
where, $\nabla$
is the Cartan h-covariant derivative,
 then  depending on the number of critical points of $\rho$, i.e. zero, one or two respectively,  it is conformal to\\
 \textbf{(a)}   A direct product $J\times\overline{M}$ of an open
interval $J$ of the real line and an       $(n-1)$-dimensional
complete
Finsler manifold $\overline{M}$.\\
 \textbf{(b)} An n-dimensional Euclidean space.\\
 \textbf{(c)} An n-dimensional  unit sphere in an Euclidean
 space.}

Here, we show that if $(M,g)$ is  compact then only the third case may occur.
More precisely we prove;\\

\emph{\textbf{Theorem 1}}: \emph{Let $(M,g)$  be an  $n$-dimensional ($n\geq 2$) without boundary compact simply connected
Finsler manifold.
 In order that $(M,g)$ admits  a non-trivial solution $\rho$ of Eq.(\ref{c-field}), it is
 necessary and sufficient that  $(M,g)$ be conformally homeomorphic to  an standard n-sphere in
 the Euclidean space $\R^{n+1}$.}\\

\emph{\textbf{Theorem 2}}: \emph{Let $(M,g)$ be an $n$-dimensional ($n\geq 2$) without boundary compact simply connected
 Finsler manifold  of positive constant flag curvature. Then it is  conformally homeomorphic to  an  standard n-sphere in the Euclidean space $\R^{n+1}$.}
\section{Preliminaries.}  Let $M$ be a real n-dimensional  manifold of class $C^ \infty$. We
denote by  $TM\rightarrow M$ the  bundle
  of tangent vectors and by $ \pi:TM_{0}\rightarrow M$ the fiber bundle of
non-zero tangent vectors.
  A {\it{Finsler structure}} on $M$ is a function
$F:TM \rightarrow [0,\infty )$, with the following properties:
$F$ is differentiable ($C^ \infty$) on $TM_{0}$;  $F$ is
positively homogeneous of degree one in $y$, i.e.
 $F(x,\lambda y)=\lambda F(x,y),  \forall\lambda>0$, where we denote
 an element of $TM$ by $(x,y)$.
 The Hessian matrix of $F^{2}$ is positive definite on
$TM_{0}$, i.e.  $(g_{ij}):=\left({1 \over 2} \left[
\frac{\partial^{2}}{\partial y^{i}\partial y^{j}} F^2
\right]\right).$
 A \textit{Finsler
manifold} $(M,g)$ is a pair of a differentiable manifold $M$ and a  tensor field $g=(g_{ij})$.
Let $(x,y)$ be the line element of $TM$ and $P(y,X)\subset T_{x}(M)$  a 2-plane
  generated by the vectors $y$ and $X$ in
$T_{x}M$. Then the \emph{sectional } or \emph{flag} curvature $K(x,y,X)$ with respect
to  the plane $P(y,X)$ at a point $x\in M$ is defined by
$K(x,y,X):=\frac{g(R(X,y)y,X)}{g(X,X)g(y,y)-g(X,y)^{2}},$
where $R(X,y)y$ is the $h$-curvature tensor of Cartan connection. If
$K$ is independent of $X$, then $(M,g)$ is called \emph{space of scalar curvature}. If $K$ has no dependence on $x$ or $y$,
then the Finsler manifold is said to be of \emph{constant curvature}.
We say that a curve $\gamma$ on $M$ is a \textit{geodesic} of a Finsler
connection $\nabla$, if  its natural lift $\tilde \gamma$ to $TM$,  is a
geodesic of $\nabla$, or equivalently $\nabla_{\frac{d\tilde
\gamma}{dt}}\frac{d\tilde \gamma}{dt} =0. $
%Let $p$ and $q$ be two points on a geodesic  $\gamma$  on the Finsler manifold $(M,g)$.
%If $\gamma$ is minimal from $p$ to a certain point $q$ but not
%beyond, then $q$ is called the \emph{cut point} of $p$ along the said geodesic.
Two points $p$ and $q$ are said to be \emph{conjugate points} along a geodesic $\gamma$ if
 there exists a non-zero Jacobi field along $\gamma$ that vanishes at $p$ and $q$, cf. \cite{BCS}.
Throughout this paper, all manifolds are supposed to be connected.
%\subsection{Adapted coordinates}
 Let $\rho:M\rightarrow [0,\infty)$  be a scalar function on $M$ and
%  \begin{equation}
% \label{c-field}
$\nabla_i\nabla_j\rho= \phi g_{ij},$
% \end{equation}
a  second order differential equation, where $\nabla_i$ is the Cartan horizontal covariant derivative  and
$\phi$ is a function of $x$ alone, then we say that Eq. (\ref{c-field}) has a solution $\rho$. The solution $\rho$ is said to be  \emph{trivial} if it is constant.
Existence of solution of
 Eq. (\ref{c-field}) is equivalent to the existence of some special conformal change of metric on $M$.
We denote by $\verb"grad" \rho=\rho^{i}{\partial}/{\partial
x^i}$ the gradient
vector field of $ \rho$, where  $\rho^i = g^{ij} \rho_j$, $\rho_j ={\partial
\rho}/{\partial x^j}$ and $i, j, ...$ run over the range $1,...,n$.
We say the point $o$ of $(M,g)$ is a \emph{critical point} of $\rho$ if the
vector field $\verb"grad" \rho$ vanishes at  $o$, or equivalently if $
\rho'(o)=0$, where  $\rho'=d\rho/dt$. All other points are called
  \emph{ordinary points} of  $\rho$ on $M$.
It's noteworthy to recall that the partial derivatives $\rho_j $ are defined on the
 manifold $M$, while $ \rho^{i}$ the components of $\verb"grad" \rho$ are defined
on the slit tangent bundle $TM_0$. Hence, $\verb"grad" \rho$
 can be considered as a section of $\pi^*TM\rightarrow TM_0$,
the pulled-back tangent bundle over $TM_0$, and  its trajectories lie on $TM_0$.
 Let the Finsler manifold $(M,g)$  admits a non-trivial solution $\rho$ of (\ref{c-field}),
 then for any ordinary point $p\in M$
 there exists a coordinate neighborhood $ \mathcal{U}$ of  $p$ which contains no critical
point, and where we can choose a system of coordinates $(u^1=t, u^2,...,u^n)$ having the
 following properties, cf. \cite{AB};

 - The function $\rho$ depends only on the first variable $u^1=t$ on $\mathcal{U}$.

 - The integral curve of $\verb"grad" \rho$ is a geodesic and geodesic containing such a curve is
 called a \emph{$\rho$-curve} or a \emph{$t$-geodesic} of $\rho$.

 - The connected component of a regular hyper-surface defined by $\rho=constant$, is called
a \emph{level set of $\rho$} or simply a $t$-level. Given a solution $\rho$ and a point
$q\in  \mathcal{U}$, there exists one and only one $t$-level set of $\rho$ passing through
 $q$. The $t$-geodesics form the normal congruence to the family of $t$-level sets of $\rho$.

 - The curves defined by $u^\alpha=$const are $t$-geodesics of $\rho$, and the parameter
 $u^1=t$ may be regarded as the arc-length parameter of $t$-geodesics.

- The components $g_{ij}$ of
  the Finsler metric tensor  $g$ satisfy
 $  g_{\alpha 1}=g_{1 \alpha}=0$, where the Greek indices $\alpha,$ run over the
range $2,3,..., n$ and the Latin indices $i,j,$ run over the
range $1,2,..., n$.

- In adapted coordinates the first fundamental form of $(M,g)$ is given by
  \begin{equation} \label{meter}
ds^{2}=(dt)^{2}+ \rho'^{2}f_{\gamma\beta}du^{\gamma}du^{\beta},
  \end{equation}
  where $f_{\gamma\beta}$ given by $g_{\gamma\beta}=\rho'^{2}f_{\gamma\beta}$ are components
 of a  metric tensor
 on a $t$-level of $\rho$ and  $g_{\gamma\beta}$ is the induced
 metric tensor of this $t$-level.
For more details about our purpose on adapted coordinates, one can refer to \cite{AB,B,T}.
\section{Compact Finsler spaces of constant curvature }
\emph{\textbf{Proof of Theorem 1.}}
Let $(M,g)$  be a an  $n$-dimensional $n\geq 2$ Finsler
manifold which admits a non trivial $C^\infty$ solution $\rho$ of Eq.(\ref{c-field}).
 Consider the so called $t$-geodesic which is integral curve of the gradient vector field $\verb"grad" \rho$ on $M$. It is well known that every $t$-geodesic is a geodesic on $M$.
 Since $M$ is compact by extension of Extreme Value Theorem to differentiable manifolds every solution $\rho$ of Eq.(\ref{c-field}) is bounded and attains its extremum  values on $M$. Once the assumption is made as $M$ is without boundary,  differentiability of  $\rho$ requires that these extremum values are  critical points.
 %It remains to check that these critical points do not coincide periodically.
 Let $O$ be a critical point for a $t$-geodesic  on $M$.
 By compactness, $M$  must have finite diameter $D$ and no $t$-geodesic longer than $D$ can remain minimizing.
  % By Hoph-Rinow Theorem, it is  geodesically complete and every geodesic  can be extended forward and backward indefinitely, cf. \cite{BCS} page 168.
%  Therefore, there exist at least one another critical point  on each
%$t$-geodesic longer than $D$, issuing from $O$. In not, there would be on $M$ only
%one critical point $O$ for each $t$-geodesic longer than $D$ emanating from $O$ which is a
%contradiction to the compactness of $M$ by above argument.
 Thus every $t$-geodesic longer than $D$ issuing from $O$ contains at least  two critical points.
Before proceeding further, we shall recall that on a Finsler manifold there is no more than two critical points of $\rho$ on every $t$-geodesic emanating from $O$, cf. \cite{AB}.
% page 439.
 Therefore, every $t$-geodesic on $(M,g)$ contains exactly two critical points.
 Thus by means of Theorem A, $(M,g)$ is conformal to an $n$-dimensional  sphere in the
 Euclidean space $\R^{n+1}$, with the first fundamental form  (\ref{meter}).
Moreover, $M$ is assumed to be simply connected and  an extension of the
Milnor theorem to Finslerian category, cf. \cite{L},
% page 172,
  implies that  $M$ is topologically homeomorphic to an $n$-sphere.

Conversely, let  $(M,g)$  be compact and conformally homeomorphic to  the n-sphere
 $S^n\subset\R^{n+1}$. The
first fundamental form of $S^n$ is given by
\begin{equation} \label{First fun. form of S^n}
g_{_{S^n}}=dt^2+ \sin^2t g_{_{S^{n-1}}},
\end{equation}
 where $g_{_{S^{n-1}}}$ is the first fundamental form
of the hypersphere $S^{n-1}$, cf. \cite{Sh2}.
Let $\gamma:= x^i(t)$ be a geodesic on $(M,g)$, by definition its differential equation is
 given by
\begin{equation} \label{geodesic eq.}
\frac{d^2x^i}{dt^2}+ \Gamma^i_{jk}\frac{dx^j}{dt}\frac{dx^k}{dt}=\varphi \frac{dx^i}{dt},
\end{equation}
where $t$ is an arbitrary parameter and $\varphi$ is a function of $t$.
If we denote
%$\sigma^i$ the vector fields on $TM$ for which projection of its integral curve
%on $M$ is $\gamma$,  then we can put
 $\frac{dx^j}{dt}=\gamma^j$, by virtue of
Eq.(\ref{geodesic eq.}) we have
\begin{equation} \label{geodesic eq.+}
\gamma^k\frac{d\gamma^l}{dx^k}+ \Gamma^l_{jk}\gamma^j\gamma^k=\varphi \gamma^l.
\end{equation}
This is equivalent to $\gamma^k(\nabla_k\gamma^l)=\varphi \gamma^l$, where
 $\nabla_k$ is the Cartan h-covariant derivative.
Denoting  $\gamma_i:=g_{il}\gamma^l$ and contracting with $g_{il}$, we obtain
$\gamma^k(\nabla_k\gamma_i)=\varphi \gamma_i,$ which leads to
\begin{equation} \label{concircular eq.}
\gamma^k(\nabla_k\gamma_i-\varphi g_{ik})=0.
\end{equation}
Conformal assumption of  $(M,g)$ to the standard sphere $(S^n,g_{_{S^n}})$ implies that the Finsler metric $g$ is positively proportional to $g_{_{S^n}}$, that is $g =e^{2\psi} g_{_{S^n}}$
 %by a positive function there is a conformal diffeomorphism $f$ for which  $f^*g =e^{\psi} g_{_{S^n}}$,
where, by the  Knebelman   theorem  $\psi$ is a  function on $M$. Therefore $g$ is also a  function on $M$ and hence a Riemannian metric.
By compactness of $M$ the vector field $\gamma^k$ is complete and Eq. (\ref{concircular eq.}) leads to $\nabla_k\gamma_i=\varphi g_{ik}$ which is equivalent to Eq. (\ref{c-field}).
This completes the proof of Theorem. \hspace{\stretch{1}}$\Box$\\
  \begin{proposition} \label{prop;complete constant curvature}
Let $(M,g)$ be an n-dimensional compact  Finsler manifold of
constant flag curvature $K$. Then the  following SODE
\begin{equation} \label{Eq;ODE}
\frac{d^2\rho}{dt}+K \rho=0,
\end{equation}
admits a non-trivial solution  on $(M,g)$  if and only if  $K>0$.
  \end{proposition}
\emph{ Proof.}
Let $(M,g)$  be a Finsler manifold of constant flag curvature $K$, then the following equation holds well, cf. \cite{AZ}, see also \cite{BCS}.
\begin{equation} \label{cfc}
\ddot{A}_{ijk}+KA_{ijk}=0,
 \end{equation}
where $A_{ijk}$ is the Cartan torsion tensor, $\dot{A}_{ijk}:=(\nabla_sA_{ijk})\ell^s$ and
$\ddot{A}_{ijk}:=(\nabla_s\nabla_tA_{ijk})\ell^s\ell^t$ and $\ell^i:=y^i/F$.
Let
 $\gamma:\R \rightarrow M$ be any geodesic parameterized by arc length  $t$ on $(M,g)$ passing through $\gamma(0)=p$, having tangent vector ${{d\gamma} \over {dt}}=\ell^i$ and  the canonical lift $\hat \gamma:={d\gamma \over{dt}}$  to $TM_0$.
%We assume that the sections $X,Y, Z \in \pi^*TM$ are fixed at $v \in I_xM$, where $I_xM=\{w \in T_xM,F(w)=1\}$ is the indicatrix.
To differentiate the Cartan torsion tensor along $\gamma$, we consider the linearly independent parallel sections  $X(t)$, $Y(t)$ and $Z(t)$ of $\pi^*TM$  along  $\hat \gamma$ with $X(0)=X$, $Y(0)=Y$ and $Z(0)=Z$ at the point $p$.
By a direct computation,  for two linearly independent parallel vector fields $X(t)$ and $Y(t)$ along $\gamma$, we have
$ \frac{d}{dt}g_{_{\dot{\alpha}(t)}}(X(t),Y(t))=0,$
see for instance \cite{AZ} or \cite{BCS}.
In this sense, a $g$-orthonormal basis for $\pi^*TM$ remains $g$-orthonormal at every point $(x(t),y(t))$ along  $\hat\gamma$.
 Therefore,   by
assuming  $A(t)=A(X(t),Y(t),Z(t))$, $\dot A(t)=\dot
A(X(t),Y(t),Z(t))$ and $\ddot A(t)=\ddot A(X(t),Y(t),Z(t))$  along $\gamma$, we have
$ \frac{d  A}{dt}=\dot A$, $ \frac{d \dot A}{dt}=\ddot A$ and  Eq. (\ref{cfc}) becomes
\begin{equation} \label{cfc2}
\frac{d^2A(t)}{dt^2}+K A(t)=0.
\end{equation}
 By assuming $\rho(t)=A(t)$ we obtain Eq. (\ref{Eq;ODE}). To complete the proof we consider three cases.\\
 \emph{Case $K>0$.}
  In this case $\rho(t)=a \cos \sqrt{K}t+b \sin \sqrt{K}t$  is a non trivial general solution for Eq. (\ref{Eq;ODE}).\\
  \emph{Case $K<0$.}
In this case  the general solution is given by
\begin{equation}\label{Eq;ODEsolution}
A(t)=\alpha e^{\sqrt{-K}t}+\beta e^{-\sqrt{-K}t}.
\end{equation}
For $v \in TM_0$, assume that the norm of Cartan torsion tensor is $\| A \|_v:=sup\ A(X,Y,Z)$,  where the supremum  is taken over all unit vectors of $ \pi_v ^* TM$. Suppose that $S_xM=\{w \in
T_xM,\  F(w)=1\}$ is the indicatrix and $\| A \|=sup_{_{v \in SM}}\|
A\|_v$,  where $SM=\bigcup_{x \in M}S_xM$.
 Since $M$ is compact the norm
$\| A\|$ is bounded. On the other hand  $M$ is compact and therefore geodesically complete and the parameter $t$  takes all the values in $(-\infty,+\infty)$. Letting $t\rightarrow +\infty$ or
$t\rightarrow -\infty$, then Eq. (\ref{Eq;ODEsolution}) implies that $A(0)=0$. In fact as $t$ approaches to $t\rightarrow \pm\infty$ the left hand side of the equation is bounded and the right hand side is infinity, so Eq. (\ref{Eq;ODEsolution}) can be hold only if the coefficients $\alpha$ and $\beta$ vanishe. Replacing it on the equation (\ref{Eq;ODEsolution}), we obtain $A(t)=0$.  That is to say the solution $\rho(t)=A(t)$ of Eq. (\ref{Eq;ODE})   is trivial.\\
\emph{Case $K=0$.}
In this case  the general solution of Eq. (\ref{Eq;ODE}) is given by
$
A(t)=\alpha +\beta t,
$
where $\alpha$ and  $\beta$ are constant.
Following the procedure described above we obtain
%By a  similar  discussions we obtain
 $\beta=0$ which implies that  $A(t)=\alpha,$ is constant and hence a trivial solution of Eq. (\ref{Eq;ODE}).
This completes proof of the proposition. \hspace{\stretch{1}}$\Box$\\ \\
\emph{ \textbf{Proof of Theorem 2.}}
Let $(M,g)$ be a  compact  Finsler manifold of  positive constant flag curvature $K$.
If we assume $\phi=-K\rho$ then Eq. (\ref{c-field}) reduces to
 \begin{equation} \label{c-field2}
\nabla_i\nabla_j\rho+K\rho g_{ij}=0.
 \end{equation}
  Following an argument similar to the one in the proof of above proposition Eq. (\ref{c-field2}) reduces to   Eq. (\ref{Eq;ODE}) along geodesics.
 By virtue of Proposition \ref{prop;complete constant curvature} there is a non-trivial solution say $\rho$ for Eq. (\ref{Eq;ODE}) and hence for  Eq. (\ref{c-field2}) on $M$.
 Therefore $(M,g)$ admits  a non-trivial solution $\rho$ for Eq.(\ref{c-field}).
  % A moment's thought shows that every non-trivial solution  $\rho$ of (\ref{Eq;ODE}),  is also a non-constant solution for (\ref{c-field2}) and hence for an special case of Eq. (\ref{c-field}).
%%On the other hand $(M,g)$ is supposed to be of constant positive flag curvature $K$  and
% By virtue of Proposition \ref{prop;complete constant curvature} there is a non-trivial solution  for the SODE (\ref{Eq;ODE}) on $M$, which is also a solution  for an special case of (\ref{c-field}).
 A simple application of Theorem 1, completes proof of Theorem 2.
\hspace{\stretch{1}}$\Box$

\end{document}